\documentclass[a4paper,12pt,leqno]{article}
\usepackage{latexsym}
\usepackage[all]{xy}

\usepackage{amssymb} 
\usepackage{amsmath} 
\usepackage{theorem}

\def\P{{\mathbf{P}}}
\def\Z{{\mathbb{Z}}}
\def\K{{\mathbb{K}}}
\def\CC{{\mathbb{C}}}
\def\R{{\mathbb{R}}}
\def\A{{\mathcal{A}}}
\def\B{{\mathcal{B}}}

\DeclareMathOperator{\rank}{rank}
\DeclareMathOperator{\codim}{codim}

\DeclareMathOperator{\Der}{Der}

\DeclareMathOperator{\pd}{pd}

\DeclareMathOperator{\POexp}{POexp}

\numberwithin{equation}{section}

\newcommand{\owari}{\hfill$\square$}

\theoremstyle{break}
\newtheorem{theorem}{Theorem}[section]
\newtheorem{prop}[theorem]{Proposition}
\newtheorem{cor}[theorem]{Corollary}
\newtheorem{lemma}[theorem]{Lemma}
\newtheorem{define}[theorem]{Definition}

\newtheorem{rem}[theorem]{Remark}

\newtheorem{example}[theorem]{Example}

\newtheorem{conj}[theorem]{Conjecture}
\newtheorem{Question}[theorem]{Question}

\title{Plus-one generated and next to free arrangements of hyperplanes}
\author{Takuro Abe
\footnote
{
Institute of Mathematics for Industry,
Kyushu University,
Fukuoka 819-0395, Japan.
Email:abe@imi.kyushu-u.ac.jp. Tel:+81-92-802-4479.
}
}
\date{\today}

\begin{document}

\maketitle

\begin{abstract}
We introduce a new class of arrangements of hyperplanes, called (strictly) plus-one generated  arrangements, from 
algebraic point of view. Plus-one generatedness is close to freeness, i.e., plus-one generated arrangements have 
their logarithmic derivation modules generated by dimension plus one elements, with relations containing one linear 
form coefficient. We show that strictly plus-one generated arrangements can be obtained if we delete a hyperplane from free arrangements. 
We show a relative freeness criterion in terms of plus-one generatedness. 
In particular, for plane arrangements, 
we show that a free arrangement is in fact surrounded by free or strictly plus-one generated arrangements. We also give several applications.
\end{abstract}

\section{Introduction}
To state the main results, let us introduce a part of notations. See \S 2 for details. 
Let $\K$ be a field, $V=\K^\ell$ and $S=\mbox{Sym}^*(V^*)$ the coordinate ring. Let $x_1,\ldots,x_\ell$ 
be a basis for $V^*$ so that $S=\K[x_1,\ldots,x_\ell]$, and that $\Der S=\bigoplus_{i=1}^\ell S \partial_{x_i}$. An 
\textbf{arrangement of hyperplanes} $\A$ in $V$ is a finite set of linear hyperplanes in $V$. We assume that every $\A$ is essential unless otherwise specified, i.e., 
$\cap_{H \in \A}H= \{0\}$. For $H \in \A$, let 
$$
\A^H:=\{L \cap H \mid L \in \A \setminus \{H\}\}
$$
be the restriction of $\A$ onto $H$, 
which is an arrangement in $H$. 
For each $H \in \A$ 
fix a linear form $\alpha_H \in V^*$ such that $\ker \alpha_H=H$. Then the \textbf{logarithmic derivation module} 
$D(\A)$ of $\A$ is defined by 
$$
D(\A):=\{\theta \in \Der S\mid \theta(\alpha_H) \in S\alpha_H\ (\forall H \in \A)\}.
$$
In general, $D(\A)$ is a reflexive graded $S$-module, and not free. So we say that $\A$ is \textbf{free} 
with \textbf{exponents} $\exp(\A)=(d_1,\ldots,d_\ell)$ if there exist homogeneous derivations $\theta_1,\ldots,\theta_\ell \in D(\A)$ such that 
$D(\A)=\bigoplus_{i=1}^\ell S\theta_i$ and that $\deg \theta_i=d_i\ (i=1,\ldots,\ell)$. Here $\theta=\sum_{i=1}^\ell f_i \partial_{x_i} \in \Der S$ is 
homogeneous of degree $d$ if all $f_i$ belongs to the homogeneous degree $d$-part $S_d$ of $S$. 

Free arrangements have been intensively studied as they are the most important algebraic properties of arrangements. However, only few of algebraic structures of 
$D(\A)$ is known when $\A$ is not free. In general, 
it is almost impossible to determine the algebraic structure of a non-free $D(\A)$ without 
using a computer program. 
In this article, we introduce a new class of arrangements from the algebraic point of view, which are next easiest to 
free arrangements as follows:

\begin{define}
%
(1)\,\,
We say that $\A$ is \textbf{plus-one generated} with \textbf{exponents} $\POexp(\A)=(d_1,\ldots,d_\ell)$ and \textbf{level} $d$ if $D(\A)$ has a minimal free resolution of the following form:
\begin{equation}
0 \rightarrow S[-d-1] \stackrel{(\alpha,f_1,\ldots,f_\ell)}{\longrightarrow} S[-d]\oplus (\oplus_{i=1}^\ell S[-d_i]) 
\rightarrow D(\A) \rightarrow 0.
\label{POGresol}
\end{equation}

(2)\,\,
We say that $\A$ is \textbf{strictly plus-one generated} with \textbf{exponents} $\POexp(\A)=(d_1,\ldots,d_\ell)$ and \textbf{level} $d$ if $\A$ is plus-one generated with the same exponents and level, and $\alpha \neq 0$ in the resolution (\ref{POGresol}).
\label{POG}
\end{define}

\begin{rem}
(1)\,\,
In other words, $\A$ is plus-one generated if there is a minimal set of homogeneous generators 
$\theta_E=\theta_1,\theta_2,\ldots,\theta_\ell$ and $\varphi$ for $D(\A)$ 
such that $\deg \theta_i=d_i,\ 
\deg \varphi=d$ and 
$$
\sum_{i=1}^\ell f_i \theta_i+\alpha \varphi=0,
$$
where $f_i \in S,\ \alpha \in V^*$. This $\alpha$ is called the \textbf{level coefficient}, and 
$\varphi$ a \textbf{level element}. If it is strict, then $\alpha \neq 0$.

(2)\,\,
When $d=d_i$ for some $i$, then the choice of the level coefficient and element is not unique. 
Hence the strict plus-one generatedness means that there is a choice of generators satisfying 
the conditions in Definiton \ref{POG} (2). 

(3)\,\,
Recall that for all $H \in \A$, there is a decomposition
\begin{equation}
D(\A) \simeq S\theta_E \oplus D_H(\A),
\label{dec}
\end{equation}
where $\theta_E$ is the Euler vector field and $D_H(\A):=\{
\theta \in D(\A) \mid \theta(\alpha_H)=0\}$. See Lemma 1.33 in \cite{Y3} for example. 
Hence choosing the generators appropriately, we may assume that the relation (\ref{POGresol}) exists in $D_H(\A)$. 

(4)\,\,
By definition, for a plus-one generated arrangement $\A$, its exponents and level cannot be computed by its characteristic polynomial. Conversely, we can compute the first and second Betti numbers. See Proposition \ref{Betti}. 
\label{remim}
\end{rem}

As the definition says, plus-one generated arrangements are not free, and have simple algebraic structures, which is very close to that of free arrangements. If $d=\max_i \{d_i\}$ and 
$\ell=3,4$, then it coincides with the definition of nearly free line and plane arrangements defined by Dimca and Sticularu in \cite{DiS1} and \cite{DiS2}. Definition \ref{POG} is motivated by them, and gives a generalization of them. 
On the other hand, we can define the arrangement close to the free arrangement in 
the sense of the inclusion relation as follows:

\begin{define}

(1)\,\, 
We say that $\A$ is \textbf{next to free minus (NT-free minus)} if 
there is a free arrangement $\B$ and $H \in \B$ such that 
$\A=\B \setminus \{H\}$.

(2)\,\, 
We say that $\A$ is \textbf{next to free plus (NT-free plus)} if 
there is a free arrangement $\B$ and a hyperplane $H \not \in \B$ such that 
$\A=\B \cup \{H\}$.

(3) $\A$ is \textbf{next to free} if $\A$ is either next to free plus or minus. 
\label{ntf}
\end{define}

Since both arrangements above are close to free arrangements in different sense,  
it is natural to ask whether 
they are related. 
In fact, we can show that next to free minus arrangements are free, or strictly plus-one generated:

\begin{theorem}
Let $\A$ be free with $\exp(\A)=(d_1,\ldots,d_\ell)$ and $H \in \A$. Then $\A':=
\A \setminus \{H\}$ is free, or strictly plus-one generated with exponents $(d_1,\ldots,d_\ell)$ and level $d=|\A'|-|\A^H|$. 
\label{del2}
\end{theorem}

Theorem \ref{del2} says that, if an arrangement is close to the free arrangement in the sense of the inclusion relation, then so is as an algebraic structure. As an example, we give an arrangement which is not free but whose algebraic structure can be
determined by Theorem \ref{del2} without any algebraic 
computation. 

\begin{example}
Let $\ell \ge 2$ and $\A$ be an arrangement in $\R^{\ell+1}$ defined by 
$$
Q(\A)=z\prod_{i=1}^\ell x_i \prod_{i=1}^\ell (x_i-z)\prod_{1 \le i < j \le \ell} (x_i^2-x_j^2)
\prod_{1 \le i < j \le \ell} (x_i-x_j-z)\prod_{1 \le i < j \le \ell} (x_i+x_j-z).
$$
This is the Shi arrangement of the type $B_\ell$, which is free with $\exp(\A)=(1,(2\ell)^\ell)$ by 
\cite{Y1}. Here 
$(a)^b$ stands for the $b$-copies of $a$ for nonnegative integers $a,b$. Let $H=\{z=0\}$ and $\A':=\A \setminus \{H\}$. Then 
$|\A^H|=\ell^2$. Thus $|\A|-|\A^H| =\ell^2+1 \not \in \exp(\A)$, and $\A'$ is not free. 
By Theorem \ref{del2}, $\A$ is strictly plus-one 
generated with $\POexp(\A')=(1,(2\ell)^\ell)$ and level $\ell^2$, i.e.,
$$
0 \rightarrow S[-\ell^2-1] \rightarrow S[-\ell^2]\oplus S[-1]
\oplus S[-2\ell]^{\oplus \ell}
\rightarrow D(\A') \rightarrow 0.
$$
\end{example}

An immediate corollary of Theorem \ref{del2} is as follows.

\begin{cor}
Let $\A$ be free, $\A \ni H$ and $\A':=\A \setminus \{H\}$. Then 

(1)\,\, the projective dimension $\pd_S D(\A') \le 1$, and 

(2)\,\,
$D(\A')$ is generated by at most $(\ell+1)$-derivations.
\label{cors}
\end{cor}

A corollary of Theorem \ref{del2} is the following criterion for freeness when the level is strictly higher than exponents, which is 
relative compared with other criterions like by K. Saito in \cite{Sa}, or Yoshinaga in \cite{Y1}, \cite{Y2}, \cite{AY}. 

\begin{theorem}[Relative freeness criterion]
Let $H \in \A,\ \A':=\A \setminus \{H\}$, and let 
$|\A'|-|\A^H| =:d$. Then $\A$ is free if 
and only if either
\begin{itemize}
\item[(1)]
$\A^H$ is free and $\chi(\A^H;t) \mid \chi(\A;t)$, or 
\item[(2)]
$\A'$ is strictly plus-one generated with level element $\varphi \in D(\A')_d$ , 
and $D(\A')_d \setminus \K\varphi \subset D(\A)_d$.
\end{itemize}
In case (1), $\exp(\A) \cap \exp(\A') = \exp(\A^H)$, and 
in case (2), $\exp(\A) = \POexp(\A')$.
\label{criterion2}
\end{theorem}

Theorem \ref{criterion2} (2) requires very strong conditions which is hard to check. However, this is the first criterion of freeness in which we use the information on the triple $(\A,\A',\A^H)$. Thus if we can improve (2), then we may apply Theorems \ref{adddel} or \ref{division} to approach 
Terao's conjecture. If the level is strictly larger than exponents, then the criterion 
becomes simpler.

\begin{theorem}[Relative freeness criterion with large level]
Let $H \in \A,\ \A':=\A \setminus \{H\}$. 
Assume that 
$\chi(\A;t)=\prod_{i=1}^\ell (t-d_i)$, and $|\A'|-|\A^H| =:d>d_i\ (\forall i)$. Then $\A$ is free if 
and only if 
$\A'$ is strictly plus-one generated with $\POexp(\A)=(d_1,\ldots,d_\ell)$.
In this case, $\exp(\A) = \POexp(\A')$ and the level is $d$.
\label{criterion}
\end{theorem}


We have the addition version of Theorem \ref{del2} as follows:

\begin{theorem}
Let $\A$ be an arrangement, $H \in \A$ and 
$\A':=\A \setminus \{H\}$. Assume that $\A'$ is free with 
$\exp(\A')=(d_1,\ldots,d_\ell)_\le$. If $d:=|\A'|-|\A^H|\ge d_{\ell-2}$, then 
$\A$ is free, or strictly plus-one generated with $\POexp(\A)=(d_1,\ldots,
d_{\ell-2},d_{\ell-1}+1,d_\ell+1)$ and 
level $d_{\ell-1}+d_\ell-|\A|+|\A^H|+1$.
\label{ell1}
\end{theorem}

Here $(d_1,\ldots,d_\ell)_\le$ stands for $d_1 \le \cdots \le 
d_\ell$. 

\begin{rem}
In fact, Theorem \ref{ell1} can be divided into two statements, i.e., $\A$ is free if $d_{\ell-1} \le 
d \le d_\ell$, and free or strictly plus-one generated when $d_{\ell-2} \le d <d_{\ell-1}$. The former result is due to 
Hoge in \cite{BCH}. See Theorem \ref{BCH} for details. 
\label{BCMrem}
\end{rem}

As we can see, Theorems \ref{del2} and \ref{ell1} are not symmetric. In fact, the addition does not work well as the deletion as in Example \ref{addnot}.
When $\ell=3$, contrary to Example \ref{addnot}, the complete symmetry exists, which 
gives a complete relative criterion for the freeness.

\begin{theorem}[Relative freeness criterion for plane arrangements]
Let $\ell=3$. Then the following are equivalent.

\begin{itemize}
\item[(1)]
$\A$ is free with $\exp(\A)=(1,a,b)$.

\item[(2)]
For any $H \in \A$, $\A':=\A \setminus \{H\}$ is either 

(i)\,\, free with $\exp(\A')=(1,a,b-1)$ and $|\A'|-|\A^H|=b-1$, or 

(ii)\,\,
(strictly) plus-one generated with $\POexp(\A')=(1,a,b)$ and level $d=|\A'|-|\A^H|$.

\item[(3)]
For some $H \in \A$, $\A':=\A \setminus \{H\}$ is either 

(i)\,\, free with $\exp(\A')=(1,a,b-1)$ and $|\A'|-|\A^H|=b-1$, or 

(ii)\,\,
(strictly) plus-one generated with $\POexp(\A')=(1,a,b)$ and level $d=|\A'|-|\A^H|$.

\item[(4)]
For any plane $L \not \in \A$, $\B:=\A \cup \{L\}$ is either 

(i)\,\, free with $\exp(\B)
=(1,a,b+1)$ and $|\B^H|=a+1$, or 

(ii)\,\,
(strictly) plus-one generated with $\POexp(\B)=(1,a+1,b+1)$ and level $d=|\B^H|-1$.

\item[(5)]
For some plane $L \not \in \A$, $\B:=\A \cup \{L\}$ is either 

(i)\,\, free with $\exp(\B)=(1,a,b+1)$ and $|\B^H|=a+1$, or 

(ii)\,\,
(strictly) plus-one generated with $\POexp(\B)=(1,a+1,b+1)$ and level $d=|\B^H|-1$.

\end{itemize}
\label{equal}
\end{theorem}

\begin{rem}
When $\ell=3$, plus-one generatedness coincides with strictly plus-one 
generatedness. See Proposition \ref{s}.
\end{rem}

As Theorem \ref{equal} implies, for plane arrangements, next to free arrangements are free or strictly plus-one generated. 
Also, free arrangements are surrounded by free or strictly plus-one generated arrangements. 

We say that $\A$ has a \textbf{free filtration} 
$$
\emptyset=\A_0 \subset \A_1 \subset \cdots \subset \A_n=\A
$$
if $|\A_i|=i$ and $\A_i$ is free for $i=0,\ldots,n$. It is shown in \cite{A4} that whether 
a free arrangement has a free filtration or not depends only on combinatorics. 
By using plus-one generatedness, we can show the following:

\begin{theorem}
Let $\ell \le 4$. Then whether $\A$ has a free filtration depends only on $L(\A)$.
\label{FF}
\end{theorem}

The organization of this article is as follows. \S2 is devoted for recalling previous results and definitions. In \S 3, we introduce two tricks which are often used in this article, 
the replacement and locally free tricks. In \S 4 we prove the main results except for Theorem \ref{equal}, which will be proved in \S 5. In \S 6 we consider whether a given non-free arrangement can admit an addition which makes it free. In \S7 we introduce 
a lot of interesting examples on free, next to free and plus-one generated arrangements. In \S 8 we consider another expression of the combinatoriality of the deletion theorem first proved in \cite{A4}.
\medskip

\noindent
\textbf{Acknowledgements}.
The author is grateful to M. Barakat, M. Cuntz, T. Hoge and H. Terao for 
the permission to put Theorem \ref{BCH} in this article.  Also, the author is grateful to 
Anca Macnic for a comment to this article. 
The author is 
partially supported by JSPS Grant-in-Aid for Scientific Research (B) JP16H03924, and 
Grant-in-Aid for Exploratory Research JP16K13744.

\section{Preliminaries}

In this section let us summarize several results and definitions. 
Start from combinatorics and topology of arrangements. 
Let 
$$
L(\A):=\{ \bigcap_{H \in \B} H \mid \B  \subset \A\}
$$ 
be the 
\textbf{intersection lattice} of $\A$ with a partial order induced from the 
reverse inclusion. Define 
$L_i(\A):=
\{X \in L(\A) \mid 
\codim_V X=i \}$. 
The \textbf{M\"{o}bius function} $\mu : L(\A) 
\rightarrow \Z$ is defined by $\mu(V)=1$, and by 
$$
\mu(X):=-\sum_{Y \in L(\A),\, X  \subsetneq Y \subset V} \mu(Y)
$$
for $X \in L(\A) \setminus \{V\}$. 
Define the \textbf{characteristic polynomial} $\chi(\A;t)$ by 
$$
\chi(\A;t):=\sum_{X \in L(\A)} \mu(X) t^{\dim X},
$$
and the \textbf{Poincar\'e polynomial} $\pi(\A;t)$ by 
$$
\pi(\A;t):=\sum_{X \in L(\A)} \mu(X) (-t)^{\codim X}.
$$
It is clear that 
$$
\pi(\A;t)=(-t)^\ell\chi(\A;-t^{-1}).
$$
For $X \in L(\A)$, the \textbf{localization} 
$\A_X$ of $\A$ at $X$ is defined by 
$$
\A_X:=\{H \in \A \mid 
X \subset H\},
$$
and 
the 
\textbf{restriction} 
$\A^X$ of $\A$ onto $X$ is defined by 
$$
\A^X:=\{H \cap X \mid
H \in 
\A \setminus \A_X\}.
$$
We agree that, for $H \in \A$ and $X \in L(\A^H)$, 
$\A^H_X$ stands for $(\A_X)^H$. It is easy to show that 
$\A_X^H=(\A^H)_X$ too. 

The characteristic polynomial $\chi(\A;t)$ 
is combinatorial but not 
easy to compute in general. The most useful 
inductive method to compute $\chi(\A;t)$ 
is so called the \textbf{deletion-restriction formula} (e.g., Corollary 2.57 in \cite{OT}) as follows:
$$
\chi(\A;t)=\chi(\A \setminus \{H\};t)-\chi(\A^H;t).
$$
Let $\chi(\A;t)=\sum_{i=0}^\ell (-1)^i b_{i}(\A) t^{\ell-i}$. 
When $\A \neq \emptyset$, it is known that $\chi(\A;t)$ is divisible by 
$t-1$. Define $\chi_0(\A;t):=\chi(\A;t)/(t-1):=\sum_{i=0}^{\ell-1} (-1)^i b_i^0(\A) t^{\ell-i-1}$. By the definition, $b_1(\A)=|\A|,\ b_1^0(A)=|\A|-1$ and 
\begin{equation}
b_2^0(\A)=b_2(\A)-(|\A|-1).
\label{b20}
\end{equation}
It is known that 
$b_i(\A)$ is the $i$-th Betti number of $V \setminus \cup_{H \in \A} H$ when 
$\K=\CC$.

Let $\theta_E$ be the \textbf{Euler 
derivation} $\theta_E=\sum_{i=1}^\ell x_i \partial_{x_i}$ which is always contained in $D(\A)$. Thus by 
(\ref{dec}), $\deg \theta_E=1$ is contained in $\exp(\A)$ if $\A \neq \emptyset$ is free. 
The localization $\A_X$ is free if $\A$ is free for any $X \in L(\A)$ 
(see \cite{OT}, Theorem 4.37). 
We say that $\A$ is \textbf{locally free} if 
$\A_X$ is free for any $X \in L(\A)$ with $X \neq \bigcap_{H \in \A} H$. $\A$ is locally free if and only if 
the sheaf $\widetilde{D_H(\A)}$ is a vector bundle on $\P^{\ell-1} =\P(V^*)=\mathbf{Proj}(S)$ for 
any (equivalently, some) $H \in \A$ (see \cite{Y3}, Proposition 1.20).

Define the \textbf{Euler restriction map} $\rho:D(\A) \rightarrow 
D(\A^H)$ of an arrangement $\A$, by taking modulo $\alpha_H$. Then we have the following:

\begin{prop}[e.g., Proposition 4.45 in 
\cite{OT}]
Let $H \in \A$. Then there is an exact sequence
$$
0 \rightarrow D(\A \setminus \{H\}) \stackrel{\cdot \alpha_H}{\rightarrow }
D(\A) \stackrel{\rho}{\rightarrow} D(\A^H).
$$
\label{exact}
\end{prop}

The following is the most important, and relates algebra, topology and combinatorics of free 
arrangements. 

\begin{theorem}[Terao's factorization, \cite{T2}]
Assume that $\A$ is free with $\exp(\A)=(d_1,\ldots,d_\ell)$. Then 
$\chi(\A;t)=\prod_{i=1}^\ell(t-d_i)$.
\label{Teraofactorization}
\end{theorem}

Let us collect several fundamental results on freeness.

\begin{prop}[\cite{T1}]
Let $H \in \A$. Then there is a polynomial $B$ of degree $|\A|-|\A^H|-1$ such that 
$\alpha_H \nmid B$, and 
$$
\theta(\alpha_H) \in (\alpha_H,B)
$$
for all $\theta \in D(\A \setminus \{H\})$. 
\label{B}
\end{prop}

\begin{theorem}[\cite{T1}, Terao's addition-deletion theorem]
Let $H \in \A$, $\A':=\A \setminus \{H\}$ and let $\A'':=\A^H$. Then two of the following 
imply the third:
\begin{itemize}
\item[(1)]
$\A$ is free with $\exp(\A)=(d_1,\ldots,d_{\ell-1},d_\ell)$.
\item[(2)]
$\A'$ is free with $\exp(\A')=(d_1,\ldots,d_{\ell-1},d_\ell-1)$.
\item[(3)]
$\A''$ is free with $\exp(\A'')=(d_1,\ldots,d_{\ell-1})$.
\end{itemize}
Moreover, all the three hold true if $\A$ and $\A'$ are both free.
\label{adddel}
\end{theorem}

\begin{theorem}[\cite{A2}, Division theorem]
Let $H \in \A$. Then $\A$ is free if $\A^H$ is free and $\pi(\A^H;t) \mid \pi(\A;t)$.
\label{division}
\end{theorem}

The following is the special case of Theorem \ref{division} 
when $\ell=3$.

\begin{prop}[\cite{A}, Theorem 1.1]
Let $\ell=3$ and $\chi(\A;t)=(t-1)(t-a)(t-b)$ with two real numbers $a \le b$. For a 
plane $L$ which may or may not be in $\A$, let $n_L:=|\{
H \cap L \neq \emptyset \mid H \in \A,\ H \neq L\}|$. Then $n_L \le a+1$ or $b+1 \le n_L$. Moreover, $\A$ and $\A'$ are both free if 
the equality holds. 
\label{between}
\end{prop}

Next let us recall the theory of multiarrangements. A pair $(\A,m)$ is a \textbf{multiarrangement} if 
$m:\A \rightarrow \Z_{\ge 1}$. Let 
$|m|:=\sum_{H \in \A} m(H)$ and $Q(\A,m):=
\prod_{H \in \A} \alpha_H^{m(H)}$. For two multiplicities 
$m$ and $k$ on $\A$, $ k \le m$ stands for  
$k(H) \le m(H)$ for all $H \in \A$. 
For $H \in \A$, let $\delta_H:\A \rightarrow \{0,1\}$ be the \textbf{characteristic multiplicity}  
defined by $\delta_H(L)=1$ if $L=H$, and $0$ otherwise. 
We may define its 
\textbf{logarithmic derivation module} 
$D(\A,m)$ as 
$$
D(\A,m):=\{\theta \in \Der S \mid \theta(\alpha_H) \in S \alpha_H^{m(H)}\ (\forall H \in \A)\}.
$$
$D(\A,m)$ is a reflexive graded $S$-module, but not free in general. Thus 
we may define the \textbf{freeness} and \textbf{exponents} of $(\A,m)$ 
in the same way as 
for $m \equiv 1$. By the reflexivity, $(\A,m)$ is free when $\ell=2$. We can also define the \textbf{localization} $(\A_X,m_X)$ of the multiarrangement $(\A,m)$ at $X \in L(\A)$ 
by $m_X:=m|_{\A_X}$. Then $(\A_X,m_X)$ is free if $(\A,m)$ is free as when $m \equiv 1$. The most fundamental criterion for freeness is the following.

\begin{theorem}[Saito's criterion, \cite{Sa}, \cite{Z}]
Let $\theta_1,\ldots,\theta_\ell$ be homogeneous elements in $D(\A,m)$. Then 
$D(\A,m)$ has a free basis $\theta_1,\ldots,\theta_\ell$ if and only if 
they are $S$-independent, and 
$|m|=\sum_{i=1}^\ell \deg \theta_i$. Also, for the matrix 
$M:=(\theta_i(x_j))$, the above two statements are equivalent to 
$$
\det M=cQ(\A,m)
$$
for some $c \in \K \setminus \{0\}$.
\label{Saito}
\end{theorem}

We can construct the multiarrangement canonically from an arrangement $\A$ in the 
following manner:

\begin{define}[\cite{Z}]
For an arrangement $\A$ in $\K^\ell$ and $H \in \A$, define the \textbf{Ziegler 
multiplicity} $m^H:\A^H \rightarrow \Z_{>0}$ by  
$m^H(X):=|\{L \in \A \setminus \{H\} \mid 
L \cap H=X\}|$ for $X \in \A^H$. The pair $(\A^H,m^H)$ is 
called the \textbf{Ziegler restriction} of $\A$ onto $H$. Also, there is a  
\textbf{Ziegler restriction map}
$$
\pi:D_H(\A) \rightarrow D(\A^H,m^H)
$$
by taking modulo $\alpha_H$. In particular, there is an exact sequence 
$$
0 \rightarrow D_H(\A) \stackrel{\cdot \alpha_H}{\rightarrow} 
D_H(\A) \stackrel{\pi}{\rightarrow}D(\A^H,m^H).
$$
\label{Zieglerrest}
\end{define}

By definition, $|m^H|=|\A|-1$. A remarkable property of Ziegler restriction maps is the following.

\begin{theorem}[\cite{Z}]
Assume that $\A$ is free with $\exp(\A)=(1,d_2,\ldots,d_\ell)$. Then 
for any $H \in \A$, the Ziegler restriction $(\A^H,m^H)$ is also free with 
$\exp(\A^H,m^H)=(d_2,\ldots,d_\ell)$. Explicitly, for the Ziegler restriction $\pi:D_H(\A) \rightarrow D(\A^H,m^H)$, 
there is a basis $\theta_2,\ldots,\theta_\ell$ for $D_H(\A)$ such that 
$\pi(\theta_2),\ldots,\pi(\theta_\ell)$ form a basis for $D(\A^H,m^H)$ (recall 
the decomposition (\ref{dec})). 
In particular, 
$\pi$ is surjective when $\A$ is free.
\label{Ziegler2}
\end{theorem}

The converse of Theorem \ref{Ziegler2} is proved by Yoshinaga. 
In this article we mainly use it when $\ell=3$.

\begin{theorem}[Yoshinaga's criterion, \cite{Y2}]
Let $\ell=3$ and $H \in \A$. Let 
$\chi_0(\A;t)=t^2-(|\A|-1)t+b_2^0(\A)$, and 
$\exp(\A^H,m^H)=(d_1,d_2)$. Then 
$b_2^0(\A) \ge d_1d_2$, and $\A$ is free if and only if $b_2^0(\A)=d_1d_2$. Moreover, 
$b_2^0(\A)-d_1 d_2=\dim_\K\mbox{coker}(\pi:D_H(\A) \rightarrow D(\A^H,m^H) )<+\infty$.
\label{Y}
\end{theorem}

Also, the following relation between Betti numbers and Chern classe are important.

\begin{prop}[\cite{DS}, Proposition 5.18]
$b_i^0(\A)=(-1)^ic_i(\widetilde{D_0(\A)})$ for $i=0,1,2$.
\label{DS}
\end{prop}

Now we introduce the following variant of the addition theorem.

\begin{theorem}[Multiple addition theorem (MAT), \cite{ABCHT}]
Let $\A \ni H, \ \A':=\A \setminus \{H\}$. Assume that $\A'$ is free with 
$\exp(\A')=(d_1,\ldots,d_\ell)_\le$. If $d_\ell= |\A'|-|\A^H|$, then 
$\A$ is free with $\exp(\A)=(d_1,\ldots,d_{\ell-1},d_\ell+1)$.
\label{MAT}
\end{theorem}

The following is due to Hoge (unpublished, \cite{BCH}), which is 
close to but different from MAT. We give a proof due to Hoge for the completeness.

\begin{theorem}[\cite{BCH}]
Let $\A \ni H, \ \A':=\A \setminus \{H\}$. Assume that $\A'$ is free with 
$\exp(\A')=(d_1,\ldots,d_\ell)_\le$. If $d_\ell \ge |\A'|-|\A^H|=:d \ge d_{\ell-1}$, then 
$\A$ is free with $\exp(\A)=(d_1,\ldots,d_{\ell-2},d_{\ell-1}+1,d_\ell)$. 
In fact, $d_{\ell-1}<d<d_\ell$ cannot occur.
\label{BCH}
\end{theorem}

\noindent
\textbf{Proof}. When $d=d_\ell$, this is nothing but Theorem \ref{MAT}. Assume that $d_{\ell-1} \le d<d_\ell$. Let $\theta_1,\ldots,\theta_\ell$ be a 
homogeneous basis for $D(\A')$ with $\deg \theta_i=d_i$. By Proposition \ref{B}, there exists a polynomial 
$B \in S_d$ such that $\theta(\alpha_H) \in (\alpha_H,B)$ for $\theta \in D(\A')$. Let $\theta_i(\alpha_H)=a_i \alpha_H+b_i B$ for $i=1,\ldots,\ell$ by Proposition \ref{B}. 
Assume that there is some $i < \ell$ such that $\theta_i \not \in D(\A)$. We may assume that 
$i=\ell-1$. This is the same as 
$b_{\ell-1}=c \in \K \setminus \{0\}$. By replacing $\theta_{\ell-1}$ by $\theta_{\ell-1}/c$, we may 
assume that $c=1$. Thus replacing $\theta_i$ by $\theta_i - b_i \theta_{\ell-1}$, it holds that 
$\theta_i \in D(\A)$ unless $i=\ell-1$. Hence Theorem \ref{Saito} implies that $\theta_1,\ldots,\theta_{\ell-2},\alpha_H\theta_{\ell-1},\theta_\ell$ form a basis for $D(\A)$. By Theorem \ref{adddel}, this implies that $d=d_{\ell-1}$. 

Hence we may assume that $\theta_i \in D(\A)$ if $i < \ell$. Then again Theorem \ref{Saito} implies that $\theta_1,\ldots,\theta_{\ell-1},\alpha_H\theta_\ell$ form a basis for $D(\A)$, and 
$d=d_\ell$ by Theorem \ref{adddel}, which is a contradiction. \owari


\begin{rem}
For $H \in \A$ and 
$S':=S/\alpha_HS$, let $\overline{f}$ denote the image of $f \in S$ by 
the canonical surjection $S \rightarrow S'$. 
By abuse of notations, $\overline{f}$ is also denoted by 
$\rho(f)$ or $\pi(f)$ respectively, where $\rho$ and $\pi$ are 
the Euler and Ziegler restrictions of $\A$ with respect to $H$. 
\end{rem}

\section{Replacement and locally free tricks}

In this section we introduce two important techniques used for the rest of this article frequently. The first key is the following easy lemma.

\begin{lemma}
Let $H \in \A$ and $d:=|\A|-|\A^H|$. Let $\pi:D_H(\A) \rightarrow D(\A^H,m^H)$ be the Ziegler 
restriction, and $\rho:D(\A) \rightarrow D(\A^H)$ the Euler restriction. Then 

(1)\,\,
$\rho|_{D_H(\A)}=\pi$.

(2)\,\, Let 
$t:\A^H \rightarrow \A \setminus \{H\}$
be a section such that 
$t(X) \cap H=X$ for all $X \in \A^H$. 
Define $Q':=\prod_{X \in \A^H} \alpha_{t(X)}^{m^H(X)-1} \in S_{d-1}$ and $\theta_E^H:=
\rho(Q'\theta_E)$. Then $\theta_E^H \in D(\A^H,m^H)_d$.
\label{A}
\end{lemma}

\noindent
\textbf{Proof}. Immediate by definitions. \owari
\medskip

The following has been well-known and used by specialists. Here we put it since we also use this 
frequently. 

\begin{prop}
Let $ H \in \A$, $\theta \in D(\A \setminus \{H\})$ with $\deg \theta=
d:=|\A'|-|\A^H|$. Assume that $\theta \not \in D(\A)$. Then for all $\varphi \in D(\A')$, there 
is $f \in S$ such that $\varphi - f\theta \in D(\A)$.
\label{Bapply}
\end{prop}

\noindent
\textbf{Proof}. 
Let $B \in S_d$ be the polynomial in Proposition \ref{B}. Thus 
we may put 
$\theta(\alpha_H)=g \alpha_H+B,\ \varphi(\alpha_H)=h\alpha_H+fB$ with 
$f,g,h \in S$. Then the statement is clear.\owari
\medskip

The following is a corollary of Proposition \ref{Bapply}. 

\begin{cor}
Let $H \in \A$. Assume that there exists $\varphi \in D(\A \setminus \{H\})_{|\A|-|\A^H|-1}$ such that $\varphi \not \in D(\A)$. Then 
$$
D(\A')=D(\A)+S\varphi.
$$
\label{+1}
\end{cor}

The first technique on freeness is the following \textbf{replacement trick}, which was first 
shown in \cite{AT}, Propositon 2.7. 

\begin{theorem}[Proposition 2.7, \cite{AT}]
Let $\A$ be free with $\exp(\A)=(1,d_2,\ldots,d_\ell)_\le$, $H \in \A$, and let 
$\theta_E,\theta_2,\ldots,\theta_\ell$ be a basis for $D(\A)$ 
with $\theta_i \in D_H(\A)$ and $\deg \theta_i=d_i$. 
Assume that $\theta_E^H$ can be expressed as 
$$
\theta_E^H=\sum_{i=2}^\ell \pi(f_i \theta_i),
$$
and $\pi(f_k) \neq 0$ for some $k$ with $d_k=d:=|\A|-|\A^H|$. Then $\A':=\A\setminus \{H\}$ and 
$\A^H$ are both free. 
\label{replacement}
\end{theorem}


The next one works when we consider $D(\A)$ generated by $(\ell+1)$-elements and 
$\A$ is locally free. We call this the \textbf{locally free trick}, which was first used in \cite{A4} implicitly.

\begin{theorem}
Assume that $\A$ is locally free but not free, and $D(\A)$ has a minimal free resolution of the following form:
$$
0 \rightarrow S[-d-1] \stackrel{F}{\rightarrow}
\oplus_{i=1}^{\ell+1} S[-d_i] \stackrel{G}{\rightarrow} D(\A) \rightarrow 0.
$$
Let $F=(f_1,\ldots,f_\ell,f_{\ell+1})$ and $G=(\theta_E,\theta_2,\ldots,\theta_\ell,\theta_{\ell+1})$. 
Then $f_i \neq 0$ for all $i =2,\ldots,\ell+1$.
\label{locallyfree}
\end{theorem}

\noindent
\textbf{Proof}.
Assume that $f_{\ell+1}=0$. Then the relation among $\theta_E,\theta_2,\ldots,\theta_{\ell+1}$ is of the form 
$$
f_1 \theta_E+\sum_{i=2}^\ell f_i \theta_i=0.
$$
Since the resolution is minimal, we may assume that $\deg f_i >0$ for all $i$ if $f_i 
\neq 0$. Recall that the freeness is independent of the extension of the base field $\K$ 
(see \cite{OT}, the first paragraph in page 151). So we may assume that $\K$ is an infinite field. Now consider the zero locus 
$$
Z:=\cap_{i=2}^\ell \{f_i=0\} \subset \P(V)=\P^{\ell-1}.
$$
This is not empty since it is the intersection of at most $(\ell-1)$-hypersurfaces in $\P^{\ell-1}$. Since $\theta_E$ is a nowhere vanishing vector field, it holds that $Z=Z \cap \{f_1=0\}$. Take a point $x \in Z$ and let $k_x$ be the residue field of $\mathcal{O}_{\P^{\ell-1}}$ at $x$. Note that $D(\A) \otimes k_x=k_x^\ell$ since 
$\A$ is locally free. Thus tensoring $k_x$ to the minimal free resolution, we obtain the exact sequence 
$$
k_x \stackrel{F|_x=0}{\rightarrow} k_x^{\ell+1} \rightarrow k_x^{\ell} \rightarrow 0,
$$
a contradiction. \owari

\section{Proofs of main results}

In this section we prove the main results posed in \S1. First let us 
prove fundamental properties of plus-one generated arrangements. 

\begin{prop}
Let $\A$ be plus-one generated with $\POexp(\A)=(d_1,\ldots,d_\ell)$ and level 
$d$. Then 
$b_1(\A)=|\A|=\sum_{i=2}^\ell d_i$, and $b_2(\A)=\sum_{2 \le i < j \le \ell} d_i d_j+d$.
\label{Betti}
\end{prop}

\noindent
\textbf{Proof}. Immediate from Definition \ref{POG} and Proposition \ref{DS}. \owari
\medskip

\begin{prop}
Let $\A$ be plus-one generated with the minimal set of homogeneous 
generators $\theta_1=\theta_E,\theta_2,
\ldots,\theta_\ell,\varphi=:\theta_{\ell+1}$ with the relation
$$
\sum_{i=1}^{\ell+1} f_i \theta_i=0.
$$
Let $I:=\{i\mid 1 \le i \le \ell+1,\ f_i \neq 0\}$. Then for any $j \in I$,
$$
\{\theta_i\}_{i \in I \setminus \{j\}} \cup \{\theta_k\}_{k \not \in I}
$$
are $S$-independent.
\label{indep}
\end{prop}

\noindent
\textbf{Proof}. Since there is the unique relation among $D(\A)$, and 
$\rank_SD(\A)=\ell$, it is obvious.\owari
\medskip

\noindent
\textbf{Proof of Theorem \ref{del2}}. 
Let $\theta_E=\theta_1,\theta_2,\ldots,\theta_\ell$ be a basis for 
$D(\A)$ such that $\deg \theta_i=d_i,\ 
d_1=1\le d_2 \le \cdots \le d_\ell$ and $\theta_i \in D_H(\A)\ (i \ge 2)$. Let 
$\pi:D_H(\A) \rightarrow D(\A^H,m^H)$ be the Ziegler restriction map, thus 
$\pi(\theta_2),\ldots,\pi(\theta_\ell)$ form a basis for $D(\A^H,m^H)$ by Theorem \ref{Zieglerrest}. Let 
$Q':=\prod_{X \in \A^H} \alpha_{t(X)}^{m^H(X)-1} \in S$, where 
$t:\A^H \rightarrow \A \setminus \{H\}$ is a section as in Lemma \ref{A}. Then Lemma \ref{A} implies that 
$\rho(Q'\theta_E) \in D(\A^H,m^H)_{d}$, where $d:=|\A|-|\A^H|$ and 
$\rho:D(\A) \rightarrow D(\A^H)$ is the 
Euler restriction map. Hence there are 
polynomials $f_i \in S$ such that 
$$
\rho(Q'\theta_E)=\sum_{i=2}^\ell \pi(f_i \theta_i). 
$$
Assume that $\A'$ is not free. 
By Theorem \ref{replacement}, we may assume that $d_i=0$ if $d_i=d$. Hence in the expression 
$$
\rho(Q'\theta_E)=\sum_{i=2}^\ell \pi(f_i \theta_i)=\sum_{i=2}^\ell \rho(f_i \theta_i),
$$
we have that $\pi(f_i)=0$ if $d_i \ge d$. Here we used Lemma \ref{A}. 
Let $k$ be the minimal integer such that 
$d_{k+1} \ge d$. By Proposition \ref{exact}, there is $\varphi \in D(\A')_{d-1}$ 
such that 
$$
\alpha_H \varphi=
\sum_{i=2}^k f_i \theta_i-Q'\theta_E.
$$
We prove that 
$$
D(\A')=\langle 
\theta_E,\theta_2,\ldots,\theta_\ell,\varphi\rangle_S.
$$
Note that 
$$
\sum_{i=2}^k f_i \theta_i(\alpha_H)-Q'\theta_E(\alpha_H)=-\alpha_HQ'
=\alpha_H \varphi(\alpha_H) \neq 0,
$$
and $\alpha_H \nmid Q'$. Hence $\varphi \neq 0$ and $\varphi(\alpha_H) \not \in S \alpha_H$, i.e., 
$\varphi \not \in D(\A)$. By Proposition \ref{Bapply}, 
for all $\theta \in D(\A')$, there is $f \in S$ such that 
$\theta - f \varphi \in D(\A)$. 
Hence 
$\theta \in \langle 
\theta_E,\theta_2,\ldots,\theta_\ell,\varphi\rangle_S$.
Now let us determine relations among these generators. Recall that we have one relation:
$$
\alpha_H \varphi=
\sum_{i=2}^k f_i \theta_i-Q'\theta_E.
$$
Let 
$$
g'\varphi=\sum_{i=2}^\ell g_i \theta_i-h\theta_E
$$
be another relation. 
Since the right hand side is in $D(\A)$ and $\varphi \not \in D(\A)$, we may assume that 
$g'=\alpha_H g$. Then 
\begin{eqnarray*}
g(\alpha_H\varphi)&=&g(\sum_{i=2}^k f_i \theta_i-Q'\theta_E)\\
&=&\sum_{i=2}^\ell g_i \theta_i-h\theta_E.
\end{eqnarray*}
Since $\theta_E,\theta_2,\ldots,\theta_\ell$ form a basis for $D(\A)$, it holds that 
it holds that $gf_i=g_i$ and $g Q'=h$. Hence the relation is unique. \owari
\medskip




\noindent
\textbf{Proof of Theorem \ref{criterion2}}.
Almost all parts follow from Theorems \ref{del2}, \ref{adddel} and \ref{division}. 
To complete the proof, it suffices to show that the condition (2) implies the freeness of $\A$. 
Let $\theta_1=\theta_E,\theta_2,\ldots.\theta_\ell,\varphi$ be a minimal set of 
homogeneous generators for $D(\A')$ such that $\deg \theta_i=d_i,\ 
\deg\varphi=d$ and 
$$
\sum_{i=1}^\ell f_i \theta_i+\alpha \varphi=0.
$$
By the assumption, 
$\alpha \neq 0,\ 
\theta_i \in D(\A)$ and $\varphi \not \in D(\A)$. Thus $\alpha=
\alpha_H$.  
Since $\theta_1,\ldots,\theta_\ell$ are $S$-independent by Proposition \ref{indep}, and 
$\sum_{i=1}^\ell d_i=|\A|$ by Proposition \ref{Betti}, Theorem \ref{Saito} completes the proof. \owari
\medskip

\noindent
\textbf{Proof of Theorem \ref{criterion}}.
Note that $\A'$ is not free by Theorem \ref{adddel}. The ``only if'' 
part 
follows from Theorems \ref{del2} and \ref{division}. 
It suffices to show the ``if'' part. 
Let $\theta_1=\theta_E,\theta_2,\ldots,\theta_\ell$ and $\varphi$ be a generator for 
$D(\A')$ such that $\deg \theta_i=d_i < d$. Note that $\varphi$ is a level element, but 
$\deg \varphi$ is not clear. However, since $d=|\A'|-|\A^H|>d_i$ for all $i$, Proposition \ref{B} implies that $\theta_i \in D(\A)$ for all $i$. Thus $D(\A') \supsetneq D(\A)$ shows that $\deg \varphi \ge d > d_i$. By the assumption, 
$$
\sum_{i=1}^\ell f_i \theta_i+\alpha \varphi=0\ \ (0 \neq \alpha \in V^*).
$$
Since $\theta_1,\ldots,\theta_\ell$ are $S$-independent by Proposition 
\ref{indep}, and $\sum_{i=1}^\ell 
d_i=|\A|$ by 
Proposition \ref{Betti}, 
Theorem \ref{Saito} completes the proof. \owari
\medskip


\noindent
\textbf{Proof of Theorem \ref{ell1}}.
Let $\alpha_H=:x_1$ and $\theta_1=\theta_E,\theta_2,\ldots,\theta_\ell$ form a basis for $D(\A')$ with 
$\deg \theta_i=d_i$. 
Assume that $\theta_{\ell-2} \not \in D(\A)$. Then Proposition \ref{B} implies that $d=|\A'|-|\A^H|=d_{\ell-2}$, and we may assume that $\theta_{\ell-1},\theta_\ell \in D(\A)$ by Proposition \ref{Bapply}. Hence $\A$ is free with basis $\theta_E,\theta_2,\ldots,\theta_{\ell-3},\alpha_H\theta_{\ell-2},\theta_{\ell-1},\theta_\ell$. So we may assume that $\theta_{\ell-2} \in D(\A)$. If $\theta_{\ell-1} \not \in D(\A)$ and $\theta_{\ell} \in D(\A)$, then 
Theorem \ref{BCH} confirms that $\A$ is free. The same holds if 
$\theta_{\ell} \not \in D(\A)$ and $\theta_{\ell-1}  \in D(\A)$ by Theorem \ref{MAT}.

As a conclusion, the rest case is when $\theta_{j} \in D(\A)\ (\forall j \le \ell-2)$ and $\theta_i \not \in D(\A)\ (i=\ell-1,\ell)$ by Proposition \ref{B}. 
So again by Proposition \ref{B}, we may assume that 
$d_{\ell-2} \le d < d_{\ell-1}$. Let $B$ be a polynomial of degree $d$ satisfying 
$$
D(\A')(x_1)  \subset (x_1,B)
$$
as in Proposition \ref{B}. We may assume that $B \in \K[x_2,\ldots,x_\ell]$. 
From now on, $\equiv$ implies to take modulo $x_1$. 
Let $D_1:=\oplus_{i=1}^{\ell-2} S \theta_i \subset D(\A)$. Then we have a canonical 
surjection $D(\A) \rightarrow D_1$ as a composition of 
$$
D(\A) \subset D(\A') \rightarrow D_1 \rightarrow 0
$$
since $D(\A')$ contains $D_1$ as a direct summand. Since this map has a canonical 
section too, it holds that 
$$
D(\A)=D_1 \oplus D_2 \subset D(\A')=D_1 \oplus D_2',
$$
where $D_2'=S \theta_{\ell-1}\oplus S \theta_\ell \supset D_2$. 
Let $\theta_i(x_1) \equiv f_i B\ (i=\ell-1,\ell)$. Since $\theta_i(x_1) \in (x_1,B)$, we may again assume that $f_i \in \K[x_2,\ldots,x_\ell]=:S' \subset S$. Let 
$(f_{\ell-1},f_\ell) \equiv g \in S'$ and $f_i=h_i g \ (h_i \in S')$. By the choice of $f_i$, $(f_{\ell-1},f_\ell)=g \in S$ too. Now let us show that 
$$
D_2'':=(x_1 \theta_{\ell-1},x_1 \theta_\ell, h_\ell \theta_{\ell-1}-h_{\ell-1}\theta_\ell)_S=D_2.
$$
By definition it is clear that $D_2'' \subset D_2$. It suffices to show that $D_2'' \supset D_2$. Let $D_2 \ni \theta=s \theta_{\ell-1} +t \theta_\ell$. If $x_1 \mid s$ and $x_1 \mid t$, then there is nothing to show. Assume that $x_1 \nmid 
s$. Since $\theta(x_1) \in S x_1$ and 
$\theta_i(x_1) \not \in S x_1\ (i=\ell-1,\ell)$, we have $x_1 \nmid t$. Take a modulo by $x_1$ of 
$\theta(x_1)=s \theta_{\ell-1}(x_1)+
t\theta_\ell(x_1)$ to obtain that 
$$
0\equiv s f_{\ell-1} B+t f_\ell B
$$
in $S/x_1 S$. Since $B \neq 0 \neq g$ in $S/x_1 S$, we have 
$$
s h_{\ell-1} \equiv -t h_\ell,
$$
implying that $s\equiv h_\ell u,\ t \equiv -h_{\ell-1}u$ for 
some $u \in S$, which shows that $\A$ is strictly plus-one generated with 
$\POexp(\A)=(d_1,\ldots,d_{\ell-2},d_{\ell-1}+1,d_\ell+1)$. 

Finally let us determine the level $c=\deg h_\ell \theta_{\ell-1}$. Let 
$\deg g=:a$. Then $\deg h_\ell=d_\ell-a-d$, where $\deg B=d$. 
Thus $c=d_{\ell-1}+d_\ell-a-d$. Now compute $b_2(\A)$ in two ways. First, by the deletion-restriction formula, 
$$
b_2(\A)=b_2(\A')+|\A^H|.
$$
Thus 
$$
b_2(\A)=\sum_{2 \le i< j \le \ell} d_i d_j+2\sum_{i=2}^\ell d_i+1-d.
$$
On the other hand, by Proposition \ref{Betti}, 
\begin{eqnarray*}
b_2(\A)&=&\sum_{2 \le i< j \le \ell} d_i d_j+2\sum_{i=2}^{\ell-2} d_i+d_{\ell-1}+d_\ell+1
+c\\
&=&
\sum_{2 \le i< j \le \ell} d_i d_j+2\sum_{i=2}^{\ell} d_i+1-a-d.
\end{eqnarray*}
Hence $a=0 \iff g=1 \iff (f_{\ell-1},f_\ell)=1$, and 
$$
c=d_{\ell-1}+d_\ell-d=d_{\ell-1}+d_\ell-|\A|+|\A^H|+1
$$
as desired. \owari
\owari
\medskip

An easy conclusion is the following.

\begin{cor}
Assume that $\A$ is free with $\exp(\A)=(d_1,d_2,\ldots,d_\ell)_\le$ with $d_1=1$, and 
$H \in \A$. If $\A':=\A \setminus \{H\}$ is locally free, then 
$d:=|\A|-|\A^H| =d_i $ for some $i$, or $d > d_\ell$.
\label{roots}
\end{cor}

\noindent
\textbf{Proof}. Since $\theta_E^H \in D(\A^H,m^H)_d$, it holds that $d\ge d_2$. 
Assume that $d_k<d < d_{k+1}$ for some $2 \le k \le \ell-1$. By Theorem \ref{del2}, 
$D(\A')$ is generated by $\theta_E,\theta_2,\ldots,\theta_\ell$, a basis for $D(\A)$, and 
$\varphi \in D(\A')_{d-1}$. There are the unique relation 
$$
\alpha_H \varphi=Q'\theta_E+\sum_{i=2}^\ell f_i \theta_i.
$$
By the assumption on $d=\deg \alpha_H \varphi$, at least $f_\ell=0$, which contradicts Theorem \ref{locallyfree}. \owari
\medskip

As far as we investigated, the conclusion in Corollary \ref{roots} holds true without 
the assumption of local freeness of $\A'$ (see \S7 for example). So let us pose the following conjecture.

\begin{conj}
Let $\A$ be free with $\exp(\A)=(d_1,d_2,\ldots,d_\ell)_\le$ with $d_1=1$, and 
$H \in \A$. Then 
$d:=|\A|-|\A^H| =d_i $ for some $i$, or $d > d_\ell$.
\label{largeconj}
\end{conj}

By Proposition \ref{between}, Conjecture \ref{largeconj} is true if $\ell =3$. 
We can relate the number of small roots of the characteristic polynomials with a level 
of a plus-one generated arrangement in terms of 
the local freeness at some codimension.

\begin{theorem}
Let $\A$ be free, $H \in \A$, $\A':=\A \setminus \{H\}$ and 
let $d:=|\A|-|\A^H|$. Assume that $\A'$ is locally free at codimension 
$k$ along $H$, i.e., $\A_X$ is free for all $X \in L_k(\A)$ with $H \supset X$. Then $\chi(\A;t)$ has $d$ as its root, or 
it has at least $(k+1)$-roots less than $d$.
\label{smallroots}
\end{theorem}

\noindent
\textbf{Proof}. Assume that $\chi(\A;d) \neq 0$. Let $d_1=1 \le d_2 \le 
\cdots \le d_\ell$ be roots of $\chi(\A;t)$ with $d_s<d < d_{s+1}$. 
Assume that $s < k+1$. By Theorem \ref{del2}, there are a minimal set of 
homogeneous generators $\theta_1=\theta_E,\theta_2,\ldots,\theta_\ell,\varphi$ for $D(\A')$ such that $\deg \theta_i=d_i,\ \deg \varphi=d-1$ and 
$$
f_1 \theta_E=\sum_{i=2}^s f_i \theta_i+\alpha_H \varphi
$$
at degree $d$ for $f_i \in S$ by the reason of degrees. Then 
by the same proof of Theorem \ref{locallyfree}, the zero locus $\alpha_H=f_2=\cdots=f_s=0$ is not empty in 
$\P(V)\simeq \P^{\ell-1}$, and it coincides with $f_1=\alpha_H=f_2=\cdots=f_s=0$. Take a prime ideal $P$ from that locus. By definition, 
$\codim P\le s \le k$. Let 
$$
0\rightarrow S[-d-1] \rightarrow \oplus_{i=1}^\ell S[-d_i] \oplus S[-d] 
\rightarrow D(\A') \rightarrow 0
$$
be the minimal free resolution of the plus-one generated arrangement $\A'$. 
Tensoring the residue field $k_P$ at $P$ gives us the exact sequence 
$$
k_P \stackrel{0}{\rightarrow} k^{\ell+1}_P \rightarrow k_P^\ell \rightarrow 0,
$$
a contradiction. 
\owari

\medskip

Based on these main results and Corollary \ref{cors}, it seems natural to pose the following conjectures and problems.

\begin{Question}
Let $\A \ni H$. Then

(1)\,\,  are there any relation between $\pd_SD(\A)$ and $\pd_SD(\A')$?

(2)\,\, Are there any relation between the cardinalities of minimal sets of homogeneous generators for $D(\A)$ and $D(\A')$?
\end{Question}

Finally, let us show Theorem \ref{FF}.
\medskip

\noindent
\textbf{Proof of Theorem \ref{FF}}. 
By Proposition \ref{between}, 
this is true when $\ell \le 3$. 
Hence we may 
assume that $\ell=4$. By Proposition \ref{between} again, whether $\A$ has a filtration $\{\A_i\}_{i=0}^n$ such that $\A_i$ is all locally free or not depends only on 
$L(\A)$. If there are no such locally free filtrations of $\A$, then 
$\A$ does not have a free filtration. 
Assume that there is a locally free filtration. 
Then it suffices to show that, 
if $\B:=\A\cup \{H\}$ is locally free, and $\A$ is free, then whether 
$\B$ is free or not depends only on $L(\B)$. Let $\chi_0(\A;t)=(t-a)(t-b)(t-c)$ with 
$a\le b \le c$. By the proof of Theorem \ref{ell1}, Theorems \ref{MAT} and \ref{BCH}, 
$\B$ is free if $|\A|-|\A^H| \in \{b,c\}$. Also, 
if both $\A$ and $\B$ are free, then $|\A|-|\A^H|$ is a root of $\chi(\A;t)$ by 
Theorem \ref{adddel}. Thus we may assume that 
$|\A|-|\A^H|=a<b$. Let $\theta_a,\theta_b,\theta_c$ be a basis for $D_L(\A)$ of corresponding degrees, where $L \in \A$. If $\theta_a(\alpha_H) \not \in S\alpha_H$, then by Proposition \ref{B}, we may put 
$$
\theta_a(\alpha_H)=f_a \alpha_H+B,\ 
\theta_b(\alpha_H)=f_b \alpha_H+g_b B,\ 
\theta_c(\alpha_H)=f_c \alpha_H+g_c B,
$$
where $B \in S_{a}$. Thus replacing $\theta_b,\theta_c$ by 
$\theta_b-g_b \theta_a,\theta_c-g_c \theta_a$,
we may assume that $\theta_b,\theta_c \in D(\B)$ and thus $\B$ is free with basis 
$\theta_E,\alpha_H\theta_a,\theta_b,\theta_c$, which contradicts Theorem \ref{adddel} and 
$|\A|-|\A^H|=a<b \le c$. Assume that 
$\theta_a \in D(\B)$. If $\theta_b \in D(\B)$, then 
the same argument as the above shows that $D(\B)$ is free with basis 
$\theta_E,\theta_a,\theta_b,\alpha_H\theta_c$.
Since the same hold if $\theta_c \in D(\B)$, we may 
assume that both are not in $D(\B)$. Thus we have a splitting surjection $D_L(\B) \rightarrow S \theta_a$ as a composition of 
$$
D_L(\B) \subset D_L(\A) \rightarrow S \theta_a.
$$
Hence $D_L(\B)=S\theta_a \oplus D_{b,c}$ such that 
$D_{b,c} \subset S\theta_b \oplus S \theta_c$. Since $\B$ is locally free, so is $D_{b,c}$. Now by Theorem \ref{ell1},
$D_{b,c}$ has 
a minimal free resolution:
$$
0 \rightarrow S[-d-1]\rightarrow 
S[-b-1]\oplus S[-c-1] \oplus S[-d] \rightarrow D_{b,c} \rightarrow 0.
$$
Here $b \le c \le d$ and $d$ is the level of $\B$, which contradicts Theorem \ref{locallyfree} and the local freeness of $\B$. So $\theta_a \in D_L(\B)$ cannot 
occur. \owari
\medskip

As a corollary of the proof above, we have the following.

\begin{cor}
Let $H \in \A$, 
$\A':=\A \setminus \{H\}$ and assume that $\A'$ is 
free with $\exp(\A')=(1,d_2,\ldots,d_\ell)_\le$.  
If $\A$ is locally free, then $\A$ is free if and only if 
$|\A|-|\A^H| \in \{d_{\ell-2},d_{\ell-1},d_\ell\}$.
\label{cor1}
\end{cor}

\section{Three dimensional case}

In this section, let $\ell=3$, and we prove Theorem \ref{equal} by dividing the statements and proving each of them. First note the following special facts when $\ell=3$.

\begin{prop}
Let $\ell=3$. Then every plus-one generated arrangement is 
strictly plus-one generated. 
\label{s}
\end{prop}

\noindent
\textbf{Proof}. 
Assume that $\A$ is plus-one generated with $\POexp(\A)=(1,d_1,d_2)_\le$ and level $d$, i.e., 
there is a minimal set of homogeneous generators $\theta_E,\theta_1,\theta_2,\varphi$ for 
$D(\A)$ such that $\deg \theta_i=d_i,\ \deg \varphi=d$, and 
\begin{equation}
f \theta_E+f_1 \theta_1+f_2 \theta_2+\alpha \varphi=0.
\label{eq99}
\end{equation}
Here $f,f_i \in S,\ \alpha \in V^*$. It suffices to show that $\alpha \neq 0$. 
Assume that $\alpha=0$. Subtract $f_1 \theta_1(\alpha_H)\theta_E/\alpha_H+
f_2 \theta_2(\alpha_H)\theta_E/\alpha_H$ from (\ref{eq99}) to obtain that 
$$
h \theta_E+f_1\theta_1'+f_2 \theta_2'=0,
$$
where $\theta_i':=\theta_i-
\theta_i(\alpha_H)\theta_E/\alpha_H$. Since $\theta_i'(\alpha_H)=0$, it holds that $h=0$. So we have 
$f_1 \theta_1'=-f_2 \theta_2'$. We may assume that $(f_1,f_2)=1$, which implies that $\theta:=\theta_1'/f_2 =-\theta_2'/f_1 \in D_H(\A)$. Since the resolution is minimal, 
$\deg f_i >0$ for $i=1,2$. Thus $\theta\in D_H(\A)_{<d_1}=(0)$, a contradiction. \owari
\medskip

\begin{prop}
Let $\ell=3$ and $\A$ be plus-one generated with the minimal set of 
homogeneous generators $\theta_E,\theta_1,\theta_2,\varphi$ for 
$D(\A)$ such that $\deg \theta_i=d_i,\ \deg \varphi=d$, $\theta_i,\varphi \in D_H(\A)$ for some 
$H \in \A$ and that 
$$
f_1 \theta_1+f_2 \theta_2+\alpha \varphi=0.
$$
Then $f_i \neq 0\ (i=1,2),\ \alpha \neq 0$.
\label{nonzero}
\end{prop}

\noindent
\textbf{Proof}.
By definition, $D_H(\A)$ is generated by three derivations of degrees $d_1\le d_2\le d$. 
By Proposition \ref{s}, $\alpha \neq 0$. Assume that $f_1=0$. Then by the same proof as in Proposition \ref{s}, there is $\psi \in D_H(\A)_{<d_2}$ such that 
$\theta_2,\varphi \in S_{>0} \psi$. In particular, $D_H(\A)$ is generated by $\theta_1$ and 
$\psi$, thus $\A$ is free, a contradiction. The same proof works when $f_2=0$. \owari
\medskip

\begin{prop}
Let $\ell=3$ and $\A$ be plus-one generated with $\POexp(\A)=(1,d_1,d_2)_\le$ and 
level $d$. Then $d \ge d_2$.
\label{levelmax}
\end{prop}

\noindent
\textbf{Proof}. Apply Proposition \ref{nonzero}. 
\owari
\medskip

Now let us show the deletion version as follows:

\begin{theorem}
Let $\ell=3$, $\A \ni H$ and 
$\A':=\A \setminus \{H\}$ not free.
Then $\A$ is free with exponents $(1,a,b)$ if and only if $\A'$ is strictly plus-one generated with 
$\POexp(\A')=(1,a,b)$ and level $d=|\A'|-|\A^H|$. Note that the equality $d=|\A'|-|\A^H|$ is equivalent to 
$b_2(\A)=ab+a+b$.
\label{3}
\end{theorem}

\noindent
\textbf{Proof}. 
The ``only if'' part is immediate from Theorem \ref{del2}. 
Let us prove the ``if'' part. This direction is similar to Proposition 3.4 in \cite{W}. 
Let 
$\theta_E,\theta_1,\theta_2, \varphi \in D(\A')$ generate $D(\A')$, where $\deg \theta_1=a,\ 
\deg \theta_2=b$ and $\deg \varphi=d$. 
Since $b_2(\A')=ab+d$ by Proposition \ref{Betti}, the deletion-restriction formula implies that $b_2(\A)=b_2(\A')+|\A^H|=ab+a+b$. Thus 
$\chi(\A;t)=(t-1)(t-a)(t-b)$. By Proposition \ref{between}, it holds that 
$|\A^H|\le a+1$ or $|\A^H|=b+1$, and $\A'$ is free if $|\A^H|=a+1$ or $b+1$. %
Thus we may assume that 
$|\A^H|\le a \iff d=|\A'|-|\A^H|=a+b-|\A^H| \ge b$. If $\theta_2 \not \in D(\A)$, then $D(\A') \supsetneq D(\A)$ implies that $d=b$ by Proposition \ref{B}, and 
$\varphi \in D(\A)$. Thus replacing $\theta_2$ by $\varphi$, we may assume that $\theta_E,\theta_1,\theta_2 \in D(\A)$ and $\varphi \not \in D(\A)$. Of course $\alpha_H \varphi \in D(\A)$. 
Now consider the relation 
$$
f \theta_E+f_1 \theta_1+f_2 \theta_2+\alpha \varphi=0
$$
at degree $d+1$ coming from the definition of strict plus-one generatedness and Proposition \ref{s}. Since $\theta_E,\theta_1,\theta_2 \in D(\A)$, it holds that 
$\alpha=c\alpha_H$ for $ c\in\K$. Since $\A'$ is strictly plus-one genedated, $c \neq 0$.  Thus $\alpha_H \varphi \in 
D(\A) \cap \langle \theta_E,\theta_1,\theta_2\rangle_S$. Since $\mbox{rank}_S D(\A)=3$, it holds that 
$\theta_E,\theta_1,\theta_2$ are $S$-independent by Propositions \ref{indep} and \ref{nonzero}, and the sum of their degrees coincide with $|\A|$ by Proposition \ref{Betti}. Hence Saito's criterion shows that $\A$ is free. \owari
\medskip

We can prove the addition-version.

\begin{theorem}
Let $\A$ be not free, and $H \in \A$. Then $\A':=\A \setminus \{H\}$ is free with exponents $(1,d_1,d_2)$ if and only if 
$\A$ is strictly plus-one generated with exponents $(1,d_1+1,d_2+1)$ and level $d=|\A^H|-1$.
\label{addNTfree}
\end{theorem}

\noindent
\textbf{Proof}. 
The ``only if'' part follows from Theorem \ref{ell1}. 
Let us prove the ``if'' part. 
Note that $b_2(\A)=d_1d_2+d_1+d_2+1+d$. 
Let $\theta_E, \theta_1,\theta_2,\varphi$ be the generators for $D(\A)$ such that 
$\deg \theta_1=1+d_1 \le \deg\theta_2=1+d_2 \le d=|\A^H|-1=\deg \varphi$. 
Estimate $|\A^H|=d+1$. By Proposition \ref{Betti}, $\chi_0(\A\setminus \{H\};t)=(t-d_1)(t-d_2)$. By Proposition \ref{between}, 
$1+d=|\A^H|\le d_1+1$ or $d_2+1 \le |\A^H|$. The equalities do not hold since they imply the freeness of $\A$ by Proposition \ref{between}. Hence 
we may assume that $d < d_1$ or $d_2 < d$. By Proposition \ref{levelmax}, $d_2 <d$. 
Hence $d_2+2 \le d+1=|\A^H|$. 
So $D(\A^H)$ has a free basis of degrees $1$ and $d\ (\ge d_2+1)$. Note that we have 
a degree-preserving Euler restriction map $\rho:D(\A) \rightarrow D(\A^H)$.
Now assume that $d_2+2<d+1$. Then $\rho(\theta_i)$ is of the form $ \rho(g_i\theta_E)$ for some $g_i \in S$. By replacing $\theta_i$ by $\theta_i-g_i \theta_E$ if necessary, by Proposition \ref{exact}, we may assume that $\alpha_H$ divides both $\theta_1$ and $\theta_2$.
Hence $\theta_E,\theta_1/\alpha_H,\theta_2/\alpha_H \in D(\A')$. Since they are $S$-independent 
by 
Propositions \ref{indep} and \ref{nonzero}, and $|\A|=2+d_1+d_2$, $\A'$ is free by Saito's 
criterion.

Next assume that 
$d_2+1=d$. 
First assume that $d_1+1<d_2+1=d$. Let $\psi \in D(\A^H)$ be a 
part of basis of degree $d$. If $\rho(\varphi)=\psi$, then 
$\theta_i':=\theta_i-a_i \theta_E - b_i \varphi \in \alpha_H D(\A')$ for some $a_i,b_i \in S$ by the same proof above. 
Again by Proposition \ref{nonzero}, $\theta_E,\theta_1'/\alpha_H,\theta_2'/\alpha_H$ 
is a basis for $D(\A')$. Next assume that $\rho(\varphi) = \rho(a\theta_E)$ for some 
$a \in S$. Since 
$\rho:D(\A) \rightarrow D(\A^H)$ is generically surjective, we may 
assume that $\rho(\theta_2)=\psi$. By Proposition \ref{nonzero}, we may assume that $\theta_1'$ above, and $\varphi':=\varphi-
a\theta_E$ is divisible by $\alpha_H$. Then $\theta_E,\theta_1/\alpha_H,\varphi'/\alpha_H$ 
is a basis for $D(\A')$. 

Now assume that $d_1+1=d_2+1=d$. By the above argument, replacing $\theta_1,\theta_2$ and $\varphi$ if necessary, we may assume that $\rho(\varphi)
=
\psi$. Then the same proof as above completes the proof. \owari
%
\medskip

\noindent
\textbf{Proof of Theorem \ref{equal}}. Combine Theorems \ref{3} and \ref{addNTfree}. \owari
\medskip

We may pose the following conjecture.

\begin{conj}
The (strictly) plus-one generatedness of $\A$, and its exponents and level depend only on $L(\A)$.
\label{NTFTC}
\end{conj}

We can show the relation with Terao's conjecture.

\begin{cor}
Let $\ell=3$. Then Terao's conjecture follows from Conjecture \ref{NTFTC}.
\label{TCNTFTC}
\end{cor}

\noindent
\textbf{Proof}. 
Assume that Conjecture \ref{NTFTC} is true. Assume that $\A,\ \B$ have the isomorphic intersection lattices, and $\A$ is free with $\exp(\A)=(1,d_1,d_2)$. Let $H \in \A$ and let $L \in \B$ be the corresponding hyperplane to $H \in \A$ under the lattice isomorphism. Theorem \ref{3} shows that $\A \setminus \{H\}$ is strictly plus-one generated with exponents $(1,d_1,d_2)$ and level 
$d:=|\A|-|\A^H|-1$. Then by Theorem \ref{3} and Conjecture \ref{NTFTC}, 
$\B \setminus \{L\}$ is strictly plus-one generated with the same exponents and level as 
$\A \setminus \{H\}$.
Since 
$d=|\A|-|\A^H|-1=|\B|-|\B^L|-1$, Theorem \ref{3} confirms that $\B$ is free.\owari
\medskip

\section{Possible additions to free arrangements}
First let us introduce a definition.

\begin{define}
For an arrangement $\A$, we say that $\B$ is a \textbf{free addition} of $\A$ if 
$\B \setminus \A=\{
H\}$ and $\B$ is free.
\label{freeaddition}
\end{define}

By Theorem \ref{del2}, we know that, if $\A$ admits a free addition, then $\A$ has to be 
either free or strictly plus-one generated. Then it is natural to ask how many free additions are possible. If $\A$ is free, then it is not easy to estimate. But if $\A$ is strictly plus-one generated 
with high level, then 
there is the unique free addition.

\begin{theorem}
Let $\A$ be 
next to free minus. Thus by Theorem \ref{del2}, $\A$ is strictly plus-one generated with $\POexp(\A)=(d_1,\ldots,d_\ell)$ and level $d$. 
If $d> d_i$ for all $i$, then $\A$ admits only one free addition.
\label{numberfa}
\end{theorem}

\noindent
\textbf{Proof}.
Since $\A$ is next to free minus, there is a hyperplane 
$H \not \in \A$ such that $\B:=\A \cup \{H\}$ is free. Let 
$\theta_1=\theta_E,\theta_2,\ldots,\theta_\ell$ be a basis for $D(\B)$ 
such that $\deg \theta_i=d_i$. Then by Theorem \ref{Saito}, 
for $M:=(\theta_i(x_j))$, it holds that $\det M=Q(\B)=\alpha_H Q(\A)$.
By Theorem \ref{del2}, there is $\varphi \in D(\A)_{d} \setminus D(\B)$ 
such that $\theta_1,\ldots,\theta_\ell$ together with $\varphi$ form a 
minimal set of homogeneous generators for $D(\A)$. Let $L \not \in \A$ and 
$\A_1:=\A \cup \{L\}$. Assume that $\A_1$ is free. Since $|\A_1|=|\B|=|\A|+1$ and $\deg \theta_i <\deg \varphi=d$, 
the homogeneous basis for $D(\A_1)$ can be expressed as $S$-linear 
combinations of $\theta_1,\ldots,\theta_\ell$. Since $\det M=\alpha_H Q(\A)$, it holds that 
$\B=\A_1$. \owari
\medskip

In the above proof, we used the following easy proposition.

\begin{prop}
Assume that $\A$ is strictly plus-one generated with $\POexp(\A)=(d_1,\ldots,d_\ell)$ and level $d$. 
Let $L \not \in \A$ and 
$\B:=\A \cup \{L\}$. If $d> d_i$ for all $i$, and $\B$ is free, then $\exp(\B)=(d_1,\ldots,d_\ell)$.
\label{freenumber}
\end{prop}

\noindent
\textbf{Proof}. 
Immediate by Theorem \ref{del2}. \owari
\medskip

\begin{prop}
Assume that $\A$ is strictly plus-one generated with $\POexp(\A)=(d_1,\ldots,d_\ell)$ and level $d$. 
Let $L \not \in \A$ and 
$\B:=\A \cup \{L\}$. If $d> d_i$ for all $i$, and $\B$ is free, then $|\A|-|\B^L|=d$.
\label{number2}
\end{prop}

\noindent
\textbf{Proof}. Assume that $\B$ is free. Then Proposition \ref{freenumber} implies that $\exp(\B)=(d_1,\ldots,d_\ell)$ and 
$$
b_2(\B)=\sum_{1 \le i<j \le \ell} d_i d_j,\ b_2(\A)=
\sum_{2 \le i<j \le \ell} d_i d_j+d,
$$
where $d_1=1$. Thus $|\B^L|=b_2(\B)-b_2(\A)=|\A|-d$.\owari
\medskip

\begin{Question}
Does the converse of Proposition \ref{number2} hold true? It suffices to check only when the level is contained in exponents.
\end{Question}




\section{Examples}

We collect several examples of plus-one generated.

\begin{example}
Let $\A$ be an arrangement in $\R^3$ defined by 
$$
Q(\A)=xyz(x+y+z)=0.
$$
The algebraic structure of $D(\A)$ is well-known, and let us confirm it by using Theorem \ref{equal}. Let $H:=\{x+y+z=0\}$. Since $\A':=\A \setminus \{H\}$ is free with exponents $(1,1,1)$, and $|\A^H|=3$, Theorem \ref{equal} shows that $\A$ is strictly plus-pne generated with $\POexp(\A)=(1,2,2)$ and level $2$, i.e., 
$$
0 \rightarrow S[-3] \rightarrow S[-1] \oplus S[-2]^3 \rightarrow D(\A) \rightarrow 0$$
is a minimal free resolution. This is a nearly free arrangement of planes as in \cite{DiS1}.
\label{tangent}
\end{example}

\begin{example}
Let $\A$ be defined by 
$$
Q(\A)=xyz(y-z)(x^2-y^2)(x^2-4y^2)=0.
$$
Let $H:=\{y=0\}$ and $\A':=\A \setminus \{H\}$. $\A$ is free with $\exp(\A)=(1,2,5)$. 
$\A'$ is the famous example which is not free but $\pi(\A)=(1+t)(1+3t)^2$. Since $|\A'|
-|\A^H|=5$, Theorem \ref{equal} implies that $\A'$ is strictly plus-one generated with $\POexp(\A')=(1,2,5)$ and 
level $5$, i.e., 
$$
0 \rightarrow S[-6] \rightarrow S[-1] \oplus S[-2]\oplus S[-5]^2 \rightarrow D(\A') \rightarrow 0$$
is a minimal free resolution. This is a nearly free arrangement of planes as in \cite{DiS1}.
\label{factor}
\end{example}

\begin{example}
Let $\A$ be defined by 
$$
Q(\A)=xyz(x^2-y^2)(x^2-z^2)(y^2-z^2)((y-x)^2-z^2)
((y+x)^2-z^2).
$$
Let $H:=\{y=0\}$ and $L$ be a generic line not in $\A$. $\A$ is free with $\exp(\A)=(1,5,7)$. 
Thus by Theorem \ref{equal}, $\A':=\A  \setminus \{H\}$ is 
is strictly plus-one generated with $\POexp(\A')=(1,5,7)$ and 
level $8$, i.e.,
$$
0 \rightarrow S[-9] \rightarrow S[-1] \oplus S[-5]\oplus S[-7] \oplus S[-8] \rightarrow D(\A') \rightarrow 0$$
is a minimal free resolution, so this is not nearly free. 
Also, by Theorem \ref{numberfa}, there are no plane $K \neq H$ such that 
$\A' \cup\{K\}$ is free. 

For $\B:=\A \cup \{L\}$, by Theorem \ref{equal}, 
$\B$ is 
is strictly plus-one generated with $\POexp(\A')=(1,6,8)$ and 
level $12$, i.e.,
$$
0 \rightarrow S[-13] \rightarrow S[-1] \oplus S[-6]\oplus S[-8] \oplus S[-12] \rightarrow D(\B) \rightarrow 0$$
is a minimal free resolution. This is not nearly free.
\label{B3}
\end{example}

\begin{example}
Let $\A$ be the reflection arrangement corresponding to the complex 
reflection group $G_{29}$. Then we know that, by \cite{OT}, Appendix C, 
$\A$ is free with $\exp(\A)=(1,9,13,17)$, and $\A^H$ is free with $\exp(\A^H)=(1,9,11)$ for any $H \in \A$. 
Hence $\A':=\A \setminus \{H\}$ is not free for all $H$, and 
Theorem \ref{del2} shows that $\A'$ is 
strictly plus-one generated with $\POexp(\A')=(1,9,13,17)$ and 
level $18$, i.e.,
$$
0 \rightarrow S[-19] \rightarrow S[-1] \oplus S[-9]\oplus S[-13] \oplus S[-17] \oplus S[-18]\rightarrow D(\A') \rightarrow 0$$
is a minimal free resolution. Also, by Theorem \ref{numberfa}, there are no plane $K \neq H$ such that 
$\A' \cup\{K\}$ is free. This is not nearly free.
\label{29}
\end{example}

\begin{example}
Let $\ell=4$ and $\A$ be defined as 
$$
xyzw(x-y)(y-z)(z-w)(x+w)=0
$$
in $V=\R^4$, and $H:=\{x+w=0\}$. Then $|\A|=8$ and 
$|\A^H|=6$. Moreover, $\A':=\A \setminus \{H\}$ is supersolvable, hence free with exponents 
$(1,2,2,2)$. However, it is easy to check that $\A$ is neither free nor plus-one generated, since $D(\A)$ has 
the projective dimension $2$. In fact, since $|\A'|-|\A^H|=1$, the condition in Theorem \ref{del2} is not satisfied. Moreover, $D(\A)$ has a minimal set of generators consisting of the Euler derivation, and $6$-derivations of degree $3$. 
\label{addnot}
\end{example}

\section{Application to the deletion theorem}

In \cite{A4}, the deletion theorem is proved to be combinatorial, i.e., the following is shown.

\begin{theorem}[\cite{A4}, Theorem 1.2]
Let $\A$ be free and $H \in \A$. Then $\A':=\A \setminus \{H\}$ is free if and only if 
$\pi(\A_X^H;t ) \mid \pi(\A_X;t)$ for all $X \in L(\A^H)$. Equivalently, $\A'$ is free if and only if 
for all $X \in L(\A^H)$, it holds that 
$$
b_2(\A_X)=b_2(\A_X^H)+|\A_X^H|(|\A_X|-|\A_X^H|),
$$
and $-(|\A_X|-|\A_X^H|)^{-1}$ is a root of $\pi(\A_X;t)$.
\label{A4main}
\end{theorem}

In the proof of Theorem \ref{A4main}, the algebraic structure of $D(\A^H)$ played the key role. By using the theory on the structure of $D(\A')$, we can show that Theorem \ref{A4main} can be expressed only in terms of the number of hyperplanes.

\begin{theorem}
Let $\A$ be free and $H \in \A$. Then $\A':=\A \setminus \{H\}$ is free if and only if  
$s_X:=-(|\A_X|-|\A_X^H|)^{-1}$ is a non-zero root of $\pi(\A_X;t)$ for all $X \in L(\A^H)$.
\label{A4}
\end{theorem}

\noindent
\textbf{Proof}. Induction on $\ell \ge 2$. When $\ell=2$ there is nothing to show, 
When $\ell=3$, it follows by Proposition \ref{between}. Assume that $\ell \ge 4$. Then 
by the induction hypothesis, we may assume that $\A'$ is locally free. By 
Theorem \ref{del2}, $\A'$ is either free or strictly plus-one generated, and the latter cannot occur by  
Theorem \ref{locallyfree}. \owari
\medskip

A corollary of the proof of Theorem \ref{A4} is as follows:

\begin{cor}
Let $\A$ be free and $H \in \A$. Assume that $\A'$ is 
locally free. Then $\A':=\A \setminus \{H\}$ is free if and only if  
$-(|\A|-|\A^H|)^{-1}$ is a non-zero root of $\pi(\A;t)$.
\label{A4-1}
\end{cor}

\end{document}